\documentclass[11pt]{article}
\usepackage[utf8]{inputenc}
\usepackage{amsmath,amssymb,epsfig,bbm}
\usepackage{stmaryrd,mathabx}
\usepackage{comment}
\usepackage{color}
\usepackage[T1]{fontenc}

\usepackage[textsize=small]{todonotes}
\usepackage{enumitem}
\usepackage{varwidth}
\setlist{nolistsep}
\usepackage{hyperref}

\newtheorem{defi}{Definition}[section]
\newtheorem{prop}[defi]{Proposition}
\newtheorem{theo}[defi]{Theorem}
\newtheorem{conj}[defi]{Conjecture}
\newtheorem{lemm}[defi]{Lemma}
\newtheorem{coro}[defi]{Corollary}
\newtheorem{rema}[defi]{Remark}
\newtheorem{exem}[defi]{Example}
\newtheorem{exems}[defi]{Examples}

\newcommand{\bdefi}{\begin{defi}}
\newcommand{\edefi}{\end{defi}}
\newcommand{\bprop}{\begin{prop}}
\newcommand{\eprop}{\end{prop}}
\newcommand{\btheo}{\begin{theo}}
\newcommand{\etheo}{\end{theo}}
\newcommand{\blemm}{\begin{lemm}}
\newcommand{\brema}{\begin{rema}}
\newcommand{\erema}{\end{rema}}
\newcommand{\bexer}{\begin{exem}}
\newcommand{\eexer}{\end{exem}}
\newcommand{\bexems}{\begin{exems}}
\newcommand{\eexems}{\end{exems}}
\newcommand{\bconj}{\begin{conj}}
\newcommand{\econj}{\end{conj}}
\newcommand{\elemm}{\end{lemm}}
\newcommand{\bcoro}{\begin{coro}}
\newcommand{\ecoro}{\end{coro}}
\newcommand{\dem}{\noindent{\bf Proof. }}
\newcommand{\rem}{\noindent{\bf Remark. }}

\usepackage{mathrsfs}
\renewcommand\mathcal{\mathscr}

\newcommand{\A}{{\cal A}}

\newcommand{\C}{{\cal C}}
\newcommand{\D}{{\cal D}}
\newcommand{\E}{{\cal E}}

\newcommand{\G}{{\cal G}}

\renewcommand{\L}{{\cal L}}

\newcommand{\OOO}{{\cal O}}
\renewcommand{\P}{{\cal P}}

\newcommand{\T}{{\cal T}}

\newcommand{\maths}[1]{{\mathbb #1}}  

\newcommand{\FF}{\maths{F}}
\newcommand{\HH}{\maths{H}}

\newcommand{\LL}{\maths{L}}
\newcommand{\NN}{\maths{N}}

\newcommand{\PP}{\maths{P}}

\newcommand{\RR}{\maths{R}}
\newcommand{\SSS}{\maths{S}}

\newcommand{\XX}{\maths{X}}
\newcommand{\YY}{\maths{Y}}
\newcommand{\ZZ}{\maths{Z}}

\newcommand{\ra}{\rightarrow}
\newcommand{\bs}{\backslash}

\newcommand{\wt}[1]{{\widetilde{#1}}}

\newcommand{\ga}{\gamma}
\newcommand{\Ga}{\Gamma}

\newcommand{\cqfd}{\hfill$\Box$}

\newcommand{\card}{{\operatorname{Card}}}

\newcommand{\dbs}{\backslash\!\!\backslash}

\newcommand{\Isom}{\operatorname{Isom}}

\newcommand{\Vol}{\operatorname{Vol}}
\newcommand{\vol}{\operatorname{vol}}

\newcommand{\hdr}{{\HH}^2_\RR}

\newcommand{\PSL}{\operatorname{PSL}}

\newcommand{\PGL}{\operatorname{PGL}}

\newcommand{\flow}[1]{{{\mathfrak g}^{#1}}}  

\newcounter{fig}

\title{Rate of mixing for equilibrium states \\
  in negative curvature and trees}
\author{Anne Broise-Alamichel \and Jouni Parkkonen \and
  Fr\'ed\'eric Paulin}

\begin{document}
\bibliographystyle{../alphanum}
\maketitle
\begin{abstract}
In this survey based on the book by the authors \cite{BroParPau19}, we
recall the Patterson-Sullivan construction of equilibrium states for
the geodesic flow on negatively curved orbifolds or tree quotients,
and discuss their mixing properties, emphazising the rate of mixing
for (not necessarily compact) tree quotients via coding by countable
(not necessarily finite) topological shifts. We give a new construction
of numerous nonuniform tree lattices such that the (discrete time)
geodesic flow on the tree quotient is exponentially mixing with
respect to the maximal entropy measure: we construct examples whose
tree quotients have an arbitrary space of ends or an arbitrary (at most
exponential) growth type.

\end{abstract}

\section{A Patterson-Sullivan construction of equilibrium states}
\label{sect:construction}

We refer to \cite[Chap.~3, 6, 7]{PauPolSha15} and \cite[Chap.~2, 3,
  4]{BroParPau19} for details and complements on this section.

Let $X$ be (see \cite{BroParPau19} for a more general framework)

$\bullet$~ either a complete, simply connected Riemannian manifold
$\wt M$ with dimension $m$ at least $2$ and pinched sectional
curvature at most $-1$,

$\bullet$~ or (the geometric realisation of) a simplicial tree $\XX$
whose vertex degrees are uniformly bounded and at least $3$. In this
case, we respectively denote by $E\XX$ and $V\XX$ the sets of vertices
and edges of $\XX$. For every edge $e$, we denote by $o(e),t(e),
\overline{e}$ its original vertex, terminal vertex and opposite edge.

Let us fix an indifferent basepoint $x_*$ in $\wt M$ or in $V\XX$.

Recall (see for instance \cite{BriHae99}) that a geodesic ray or line
in $X$ is an isometric map from $[0,+\infty[$ or $\RR$ respectively
into $X$, that two geodesic rays are {\it asymptotic} if they stay at
bounded distance one from the other, and that the {\it boundary at
  infinity} of $X$ is the space $\partial_\infty X$ of asymptotic
classes of geodesic rays in $X$ endowed with the quotient topology of
the compact-open topology.  When $X=\wt M$, up to a translation
factor, two asymptotic geodesic rays converge exponentially fast one
to the other, and $\partial_\infty \wt M$ is homeomorphic to the
sphere $\SSS_{m-1}$ of dimension $m-1$. When $X$ is a tree, up to a
translation factor, two asymptotic geodesic rays coincide after a
certain time, and $\partial_\infty \wt M$ is homeomorphic to a Cantor
set.

For every $x$ in $X$, the Gromov-Bourdon {\em visual distance} $d_x$
on $\partial_\infty X$ seen from $x$ (inducing the topology of
$\partial_\infty X$ ) is defined by
$$
d_x(\xi,\eta)=
\lim_{t\ra+\infty} e^{\frac{1}{2}(d(\xi_t,\,\eta_t)-d(x,\,\xi_t)-d(x,\,\eta_t))}\;,
$$
where $\xi,\eta\in\partial_\infty X$ and $t\mapsto \xi_t,\eta_t$ are
any geodesic rays converging to $\xi,\eta$ respectively. The visual
distances seen from two points of $X$ are Lipschitz equivalent.

Let $\Ga$ be a discrete group of isometries of $X$ which is {\it
  nonelementary}, that is, does not preserve a subset of cardinality
at most $2$ in $X\cup\partial_\infty X$. When $X=\wt M$, this is
equivalent to $\Ga$ being non virtually nilpotent. When $X$ is a tree,
we furthermore assume that $X$ has no nonempty proper invariant subtree
(this is not an important restriction, as one may always replace $X$
by its unique minimal nonempty invariant subtree), and that $\Ga$ does
not map an edge to its opposite one.

The {\it limit set} $\Lambda\Ga$ of $\Ga$ is the smallest nonempty
closed invariant subset of $\partial_\infty X$, which is the
complement of the orbit $\Ga x_*$ in its closure $\overline{\Ga x_*}$,
in the compactification $X\cup\partial_\infty X$ of $X$ by its
boundary at infinity.

\medskip
\noindent{\bf Examples. } (1) Let $\wt M$ be a symmetric space with
negative curvature, e.g.~the real hyperbolic plane $\hdr$, and let $\Ga$
be an arithmetic lattice in $\Isom(\wt M)$, e.g.~$\Ga = \PSL_2(\ZZ)$
acting by homographies on the upper halfplane model of $\hdr$ with
constant curvature $-1$ (see for instance \cite{Katok92}, and
\cite{Margulis91} for a huge amount of examples).

(2) For every prime power $q$, let $\XX$ be the regular tree of degree
$q+1$, and let $\Ga=\PGL_2(\FF_q[Y])$, acting on $\XX$ seen as the
Bruhat-Tits tree $\XX_q$ of $\PGL_2$ over the local field $\FF_q((Y^{-1}))$
(see for example \cite{Serre83}, and \cite{BasLub01} for a huge amount
of examples).

\begin{center}
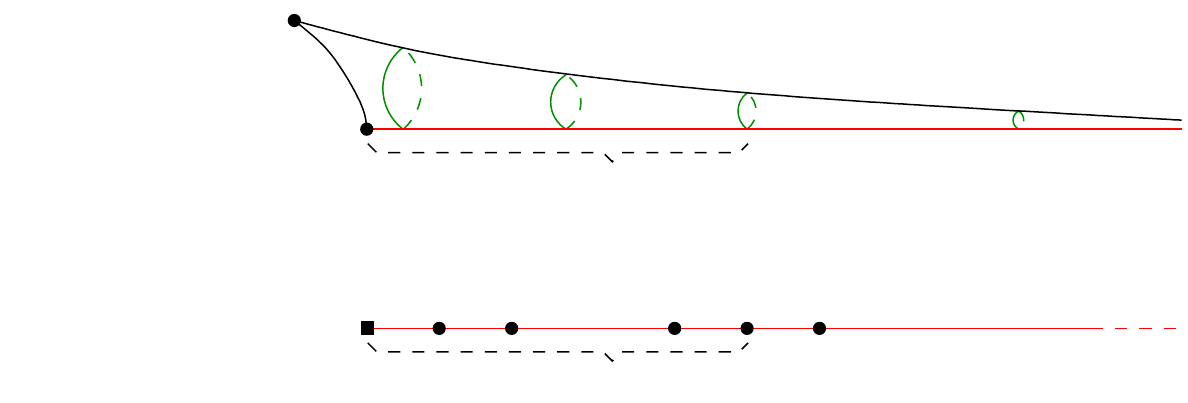
\end{center}

Note that the pictures of the quotients $\Ga\bs X$ are very similar in
the above two special examples, in particular

$\bullet$ the lengths of the closed horocycle quotients in
$\PSL_2(\ZZ)\bs \hdr$ go exponentially to $0$ (they are equal to
$e^{-t}$ where $t$ is the distance of the horocycle quotient to the
orbifold point of order $2$),

$\bullet$ the orders of the vertex stabilisers along a geodesic ray in
$\XX_q$ lifting the quotient ray $\PGL_2(\FF_q[Y])\bs\XX_q$ increase
exponentially (they are equal to $c\,q^n$ where $c$ is a constant and
$n$ is the distance of the vertex to the origin of the ray), see for
instance \cite[\S 15.2]{BroParPau19}.

\medskip
\rem Note that we allow torsion in $\Ga$, as this is in particular
important in the tree case; we allow $\Ga\bs X$ to be noncompact;
and we allow $\Ga$ not to be a lattice, which gives in the tree
case the possibility to have almost any (metrisable, compact, totally
disconnect) space of ends and almost any type of asymptotic growth of
the quotient $\Ga\bs X$ (linear, polynomial, exponential, etc), see
loc. cit.

Recall that $\Ga$ is a lattice in $X$ if either the Riemannian volume
$\Vol(\Ga\bs \wt M)$ of the quotient orbifold $\Ga\bs \wt M$ is
finite, or if the {\it graph of groups volume}
$$
\Vol(\Ga\dbs \XX)= \sum_{[x]\in\Ga\bs V\XX} \;\;\frac{1}{\card(\Ga_x)}
$$
(where $\Ga_x$ is the stabiliser of $x$ in $\Ga$) of the quotient
graph of groups $\Ga\dbs \XX$ is finite. Note the analogy, in the two
special examples above, between the computation of (most of) the
volume of $\PSL_2(\ZZ)\bs \hdr$ as a converging integral of the
lengths of the closed horocycle quotients and of the volume of
$\PGL_2(\FF_q[Y])\dbs\XX_q$ (which does converge by a geometric mean
argument).

\medskip
\noindent{\bf The phase space. }
Let $\G X$ be the space of geodesic lines $\ell:\RR\ra X$ in $X$, such
that, when $X$ is a tree, $\ell(0)$ is a vertex, endowed with the
$\Isom(X)$-invariant distance (inducing its topology) defined by
$$
d(\ell, \ell')=
\int_{-\infty}^{+\infty} d(\ell(t),\ell'(t))\;e^{-2 |t|}\,dt\;,
$$
and with the $\Isom(X)$-equivariant {\it geodesic flow}, which is
the one-parameter group of homeo\-morphisms
$$
\flow t :\ell\mapsto \{s \mapsto \ell(s+t)\}
$$
for all $\ell\in\G X$, with continuous time parameter $t\in\RR$ if
$X=\wt M$ and discrete time parameter $t\in\ZZ$ if $X$ is a tree. We
again call {\it geodesic flow} and denote by $(\flow t)_t$ the
quotient flow on the {\it phase space} $\Ga\bs \G X$.

Note that the map from the unit tangent bundle $T^1\wt M$ endowed with
Sasaki's metric to $\G \wt M$, which associates to a unit tangent vector
$v$ the unique geodesic line whose tangent vector at time $t=0$ is
$v$, is an $\Isom(\wt M)$-equivariant bi-Hölder-continuous\footnote{In
  order to deal with noncompactness issues, a map $f$ between two
  metric spaces is {\it Hölder-continuous} if there exist $c,c'>0$
  and $\alpha\in\;]0,1]$ such that for every $x,y$ in the source space,
  if $d(x,y)\leq c$, then $d(f(x),f(y))\leq c' d(x,y)^\alpha$.}
homeomorphism, by which we identify the two spaces from now on.

\medskip
\noindent
{\bf Potentials on the phase space. } We now introduce the
supplementary data (with physical origin) that we will consider on
our phase space. Assume first that $X=\wt M$. Let $\wt F:T^1\wt M\ra
\RR$ be a {\it potential}, that is, a $\Ga$-invariant,
bounded\footnote{see \cite[\S 3.2]{BroParPau19} for a weakening of
  this assumption} Hölder-continuous real map on $T^1\wt M$.
Two potentials $\wt F,\wt F^* :T^1\wt M\ra \RR$ are {\it cohomologous}
(see for instance \cite{Livsic72}) if there exists a
Hölder-continuous, bounded, differentiable along flow lines,
$\Ga$-invariant function $\wt G :T^1\wt M\ra \RR$, such that, for every
$v\in T^1\wt M$,
$$
\wt F^*(v)-\wt F(v)=\frac{d}{dt}_{\mid t=0}\wt G(\flow{t}v)\;.
$$
For every $x,y\in \wt M$, let us define (with the
obvious convention of being $0$ if $x=y$) the integral of $\wt F$
between $x$ and $y$, called the {\it amplitude} of $\wt F$ between $x$
and $y$, to be

\begin{center}
  \raisebox{0.5cm}{
$\displaystyle\int_x^y\wt F= \int_{0}^{d(x,y)} \wt F(\flow t v) \;dt$}
  ~~~~~~~~~ 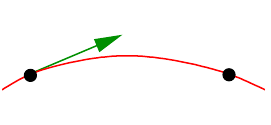
\end{center}
and $v$ is the tangent vector to the geodesic segment from $x$ to $y$.

Now assume that $X$ is a tree. Let $\wt c: E\XX\ra \RR$ be a
(logarithmic) {\it system of conductances} (see for instance
\cite{Zemanian91}), that is, a $\Ga$-invariant, bounded real map on
$E\XX$. Two systems of conductances $\wt c,\wt c^*: E\XX\ra \RR$ are
     {\it cohomologous} if there exists a $\Ga$-invariant function
     $\wt f : V\XX\ra \RR$, such that for every $e\in E\XX$
$$
\wt c^*(e)-\wt c(e)=f(t(e))-f(o(e))\;.
$$
For every $\ell\in \G X$, we denote by $e^+_0(\ell)=\ell([0,1])\in
E\XX$ the first edge followed by $\ell$, and we define $\wt F:\G X\ra
\RR$ as the map $\ell\mapsto \wt c(e^+_0(\ell))$.
For every $x,y\in V\XX$, we now define the {\it amplitude} of
$\wt F$ between $x$ and $y$, to be
\begin{center}
\raisebox{0.5cm}{
$\displaystyle\int_x^y\wt F= \sum_{i=1}^{k} \;\;\wt c(e_i) \;dt$}
  ~~~~~~~~~ 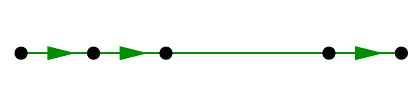
\end{center}
if $(e_1,e_2,\dots, e_k)$ is the geodesic edge path
in $\XX$ between $x$ and $y$.

In both cases, we will denote by $F :\Ga\bs \G X\ra\RR$ the function
on the phase space induced by $\wt F$ by taking the quotient modulo
$\Ga$, that we call the {\it potential} on $\Ga\bs \G X$.  Note that
we make no assumption of reversibility on $F$.

\medskip
\noindent
{\bf Cohomological invariants. } Let us now introduce three
cohomological invariants of the potentials on the phase space.

The {\it pressure} of $F$ is the physical complexity associated with
the potential $F$ defined by
\begin{center}
\fcolorbox{blue}{white}{
  $\displaystyle
  P_F= \sup_{\;\mu \;\;(\flow t)_t {\textrm{-invariant~proba~on~}}
  \Ga\bs \G X} \;\big(\; h_\mu + \int_{\Ga\bs \G X} F\;d\mu\;\big)
  $ }
\end{center}
where $h_\mu$ is the metric entropy\footnote{The metric entropy
  $h_\mu$ is the upper bound, for all measurable countable partitions
  $\xi$ of $\Ga\bs \G X$, of
  $$
  \lim_{k\ra+\infty}\;\;\frac{1}{k}\; H_\mu(\xi\vee \cdots \vee g^{-k}\xi)
  $$
  where $H_\mu(\xi)= - \sum_{E\in \xi} \mu(E) \ln \mu(E)$ is Shannon's
  entropy of the countable partition $\xi$, see for instance
  \cite{KatHas95}, and the join $\xi\vee\xi'$ of two partitions $\xi$
  and $\xi'$ is the partition by the nonempty intersections of an
  element of $\xi$ and an element of $\xi'$.} of $\mu$ for the time
$1$ map $\flow 1$ of the geodesic flow.

The {\it critical exponent} of $F$ is the weighted (by the exponential
amplitudes) orbital growth rate of the group $\Ga$, defined by
\begin{center}
  \fcolorbox{blue}{white}{
    $\displaystyle \delta_F= \lim_{n\ra+\infty}\;\frac{1}{n}\;\ln\;\Big(
    \sum_{\ga\in\Ga,\;n-1< d(x_*,\ga x_*)\leq n} \;\;
    \exp\big(\int_{x_*}^{\ga x_*} \wt F\;\big)\Big)\;. $ }
\end{center}
Note that the critical exponent $\delta_0$ of the zero potential is
the usual critical exponent of the group $\Ga$ (see for instance
\cite{Paulin97d}). We have $\delta_F\in\;]-\infty,+\infty[$ since
$$
\delta_0+\inf \wt F\leq \delta_F\leq \delta_0+\sup \wt F\;.
$$
Note that $\delta_{F\circ\iota}=\delta_F$ where $\iota : \G X\ra \G X$
is the involutive {\it time reversal map} defined by $\ell\mapsto
\{t\mapsto \ell(-t)\}$.

The {\it period} for the potential $F$ of a periodic orbit $\OOO$ of
the geodesic flow $(\flow t)_t$ on $\Ga\bs \G X$ is $\int_\OOO F=
\int_{\ell(0)}^{\ell(t_\OOO)}\;\wt F$ where $\ell\in\G X$ maps to
$\OOO$ and $t_\OOO =\inf\{t>0\;:\;\Ga \flow{t}\ell=\Ga \ell\}$ is the
{\it length} of the periodic orbit $\OOO$. The {\it Gurevich
  pressure} of $F$ is the growth rate of the exponentials of
periods for $F$ of the periodic orbits, defined by
\begin{center}
\fcolorbox{blue}{white}{ $\displaystyle \P_F^{\rm Gur}=
\lim_{n\ra+\infty}\;\frac{1}{n}\;\ln\; \sum_{\OOO\;:\; t_\OOO\leq n,\;
\OOO\cap W\neq \emptyset} \;\;\exp\big(\int_{\OOO} F\big)\;,$ }
\end{center}
where the sum is taken over the periodic orbits $\OOO$ of $(\flow
t)_t$ on $\Ga\bs \G X$ with length at most $n$ and meeting $W$, where
$W$ is any relatively compact open subset of $\Ga\bs \G X$ meeting the
nonwandering set of the geodesic flow (recall that we made no
assumption of compactness on the phase space).

Note that the above three limits exist, and are independent of the
choices of $x_*$ and $W$, and depend only on the cohomology class of
the potential $F$.

The following result proved in \cite[Theo.~4.1 and 6.1]{PauPolSha15}
extends the case of the zero potential due to Otal and Peigné
\cite{OtaPei04}.

\btheo[Paulin-Pollicott-Schapira] If $X=\wt M$ has pinched sectional
curvatures with uniformly bounded derivatives,\footnote{This
  assumption on the derivatives was forgotten in the statements of
  \cite{OtaPei04, PauPolSha15}, but is used in the proofs.} then
\begin{center}
  \fcolorbox{red}{white}{ $\displaystyle P_F=\delta_F=\P_F^{\rm Gur}$. }
\end{center}
\etheo

Note that the dynamics of the geodesic flow $(\flow t)_t$ on the phase
space $\Ga\bs \G X$ is very chaotic. In particular, there are lots
$(\flow t)_t$-invariant measures on $\Ga\bs \G X$. We give two basic
examples, and we will then contruct, using potentials, a huge family
of such measures.

\medskip
\noindent{\bf Examples. } (1) If $X=\wt M$, then the {\it Liouville
  measure} $m_{\rm Liou}$ on $T^1M=\Ga\bs (T^1\wt M)$ is the measure
on $T^1M$ which disintegrates, with respect to the canonical footpoint
projection $T^1M\ra M$, over the Riemannian measure $\vol_M$ of the
orbifold $M= \Ga\bs \wt M$, with conditional measures on the fibers
the spherical measures $\vol_{T^1_xM}$ on the (orbifold) unit tangent
spheres at the points $x$ in $M$:
\begin{center}
  \fcolorbox{blue}{white}{ $\displaystyle d m_{\rm Liou}(v)= \int_{x\in M}
     d\vol_{T^1_xM}(v)\;\;d \vol_M(x)$. }
\end{center}

(2) For every periodic orbit $\OOO$ of the geodesic flow $(\flow t)_t$
on $\Ga\bs \G X$, we denote by $\L_\OOO$ the Lebesgue
measure\footnote{If the length of $\OOO$ is $T$ and if $v\in T^1\wt M$
  maps into $\OOO$ by the canonical projection $T^1\wt M\ra T^1M$, the
  Lebesgue measure $\L_\OOO$ of $\OOO$ is the pushforward by $t\mapsto
  \Ga\flow t v$ of the Lebesgue measure on $[0,T]$.}  (when $X=\wt M$)
or counting measure (when $X$ is a tree) of $\OOO$. This is a $(\flow
t)_t$-invariant measure on $\Ga\bs \G X$ with support $\OOO$.

\medskip The main class of invariant measures we will study is the
following one, and the terminology has been mostly introduced by
Sinai, Ruelle, Bowen, see for instance \cite{Ruelle04}.  A $(\flow
t)_t$-invariant probability measure $\mu$ on the phase space $\Ga\bs
\G X$ is an {\it equilibrium state} if it realizes the upper bound
defining the pressure of $F$, that is, if
$$
h_\mu + \int_{\Ga\bs \G X} F\;d\mu = P_F\;.
$$

The remainder of this section is devoted to the problems of {\bf
  existence, uniqueness and explicit construction} of equilibrium
states.

\medskip
\noindent
{\bf Gibbs cocycles. } As for instance defined by Hamenstädt, the
(normalised) {\it Gibbs cocycle} of the potential $F$ is the function
$C:\partial_\infty X\times \wt M\times \wt M\ra \RR$ when $X=\wt M$ or
the function $C:\partial_\infty X\times V\XX\times V\XX\ra \RR$ when
$X$ is a tree, defined by the following limit of difference of
amplitudes for the renormalised potential
\begin{center}
  \raisebox{1.3cm}{
 $\displaystyle(\xi,x,y)\mapsto C_\xi(x,y)=
    \lim_{t\ra+\infty}\int_y^{\xi_t}(\wt F-\delta_F)-\int_{x}^{\xi_t}
    (\wt F -\delta_F)$,}
  ~~~~~~~~~ 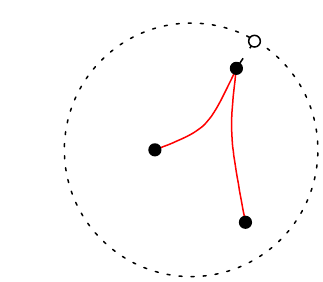
\end{center}
where $t\mapsto \xi_t$ is any geodesic ray converging to $\xi$. The
limit does exist. The Gibbs cocycle is $\Ga$-invariant (for the
diagonal action) and locally Hölder-continuous. It does satisfy the
cocycle property $C_\xi(x,z)=C_\xi(x,y)+C_\xi(y,z)$ for all
$x,y,z$. Furthermore, there exist constants $c_1, c_2>0$ (depending
only on the bounds of $\wt F$ and on the pinching of the sectional
curvature, when $X=\wt M$) such that if $d(x,y)\leq 1$, then
$C_\xi(x,y)\leq c_1d(x,y)^{c_2}$. See \cite[\S 3.4]{BroParPau19}.

\medskip
\noindent
{\bf Patterson densities. } A (normalised) {\it Patterson density} of
the potential $F$ is a $\Ga$-equiv\-ariant family
$(\mu_{x})_{x\in X}$ of pairwise absolutely continuous (positive,
Borel) measures on $\partial_\infty X$, whose support is $\Lambda\Ga$,
such that
\begin{equation}\label{eq:Pattersondensity}
\ga_*\mu_x=\mu_{\ga x}{\rm ~~~and~~~}
\frac{d\mu_x}{d\mu_y}(\xi) = e^{-C_\xi(x,\,y)}
\end{equation}
for every $\ga\in\Ga$, for all $x,y\in X$, and for
(almost) every $\xi\in\partial_\infty X$.

Patterson densities do exist and they satisfy the following Mohsen's
shadow lemma (see for instance \cite[\S 4.1]{BroParPau19}:

\smallskip\noindent
\begin{minipage}{7cm}
  \begin{center}
    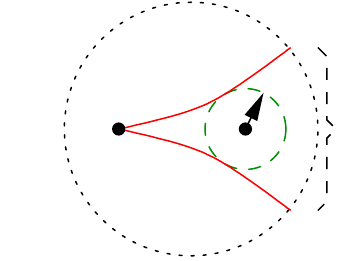~~~~~~~~~~~~~~~~
  \end{center}
\end{minipage}
\begin{minipage}{7.9cm}
Define the {\it shadow} $\OOO_xE$ seen from $x$ of a subset $E$ of $X$
as the set of points at infinity of the geodesic rays from $x$ through
$E$. Then for every $x\in X$, if $r>0$ is large enough, there exists
$\kappa>0$ such that for every $\ga\in\Ga$, we have
\end{minipage}

\begin{center}
  \raisebox{1.3cm}{\fcolorbox{green}{white}{
    $\displaystyle
    \frac{1}{\kappa}\;\exp\Big(\int_x^{\ga x}(\wt F-\delta_F)\Big)\leq
    \mu_x\big(\OOO_xB(\ga x,r)\big)\leq
    \kappa\;\exp\Big(\int_x^{\ga x}(\wt F-\delta_F)\Big)$}
    \hfill (2)}\addtocounter{equation}{1}
\end{center}

\vspace*{-0.5cm}
\noindent
{\bf Gibbs measures. } The {\it Hopf parametrisation} of $X$ at $x_*$
is the map from $\G X$ to $(\partial_\infty X\times \partial_\infty
X-{\rm Diag})\times R$, where $R=\RR$ if $X=\wt M$ and $R=\ZZ$ if $X$
is a tree, defined by
\begin{center}
  \raisebox{1.3cm}{
    \hspace{-1cm}$\displaystyle\ell\mapsto (\ell_-,\ell_+,t)$}
  ~~~~~~~~~~~~~~~~~~~~~~~~~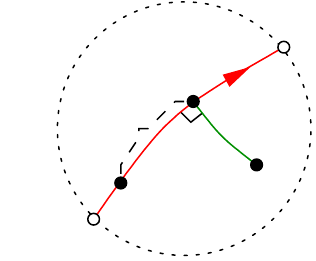
\end{center}
where $\ell_-$, $\ell_+$ are the original and terminal points at
infinity of the geodesic line $\ell$, and $t$ is the algebraic
distance along $\ell$ between the footpoint $\ell(0)$ and the closest
point to $x_*$ on the geodesic line. It is a Hölder-continuous
homeomorphism (for the previously defined distances). Up to
translations on the third factor, it does not depend on the basepoint
$x_*$ and is $\Ga$-invariant, see for instance \cite[\S 2.3 and \S
  3.1]{BroParPau19}. The geodesic flow acts by translations on the
third factor.

\medskip
Let $(\mu_{x})_{x\in X}$ and $(\mu^\iota_{x})_{x\in X}$ be Patterson
densities for the potentials $F$ and $F\circ \iota$ respectively,
where $\iota:\Ga\ell\mapsto \Ga\{t\mapsto \ell(-t)\}$ is the time
reversal on the phase space $\Ga\bs \G X$. We denote by $C^\iota$ the
Gibbs cocycle of the potential $F\circ \iota$. We denote by $dt$ the
Lebesgue or counting measure on $R$. The measure on $\G X$ defined
using the Hopf parametrisation by
\begin{center}
\fcolorbox{blue}{white}{
  $\displaystyle
  d\wt m_F(\ell)= \frac{d\mu^\iota_{x_*}(\ell_-)\;d\mu_{x_*}(\ell_+)\;dt}
  {\exp\big(\,C^\iota_{\ell_-}(x_*,\,\ell(0))+C_{\ell_+}(x_*,\,\ell(0))\,\big)}
  $ }
\end{center}
is a $\sigma$-finite nonzero measure on $\G X$. By Equation
\eqref{eq:Pattersondensity} and by the invariance of the 
measure $dt$ under translations, it is independent of the choice of
basepoint $x_*$, hence is $\Ga$-invariant and $(\flow t)_t$-invariant.
Therefore it induces a $\sigma$-finite nonzero $(\flow t)_t$-invariant
measure on $\Ga\bs\G X$, called the {\it Gibbs measure} on the phase
space and denoted by $m_F$.

\medskip
\noindent{\bf Examples. } (1) When $F=0$, then the Gibbs measure is
called the Bowen-Margulis measure (see for instance \cite{Roblin03}).

(2) When $X=\wt M$ and $\wt F$ is the {\it unstable Jacobian}, that
is, for every $v\in T^1\wt M$,
\begin{center}
  \raisebox{1.2cm}{\fcolorbox{blue}{white}{
  $\displaystyle
  \wt F^{\rm su}(v)=-\;\frac{d}{dt}_{\mid t=0}\ln\Big(
  \begin{array}{c}
  {\rm Jacobian~of~restriction~of~} \flow t{\rm ~to}\\
  {\rm strong~unstable~leaf~} W^{su}(v)
  \end{array}\Big),
  $ }}~~~~~~ 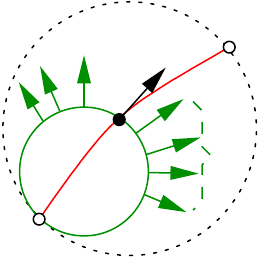
\end{center}
we have the following result (see \cite[\S 7]{PauPolSha15}, in
particular for weaker assumptions). When $M$ has variable sectional
curvature, the Liouville measure and the Bowen-Margulis measure might
be quite different. The following result in particular says that the
huge family of Gibbs measures interpolates between the Liouville
measure and the Bowen-Margulis measure. This sometimes provides common
proofs of properties satisfied by both the Liouville measure and the
Bowen-Margulis measure.

\btheo[Paulin-Pollicott-Schapira] If $X=\wt M$ has pinched sectional
curvatures with uniformly bounded derivatives, then $\wt F^{\rm su}$
is Hölder-continuous and bounded. If $\wt M$ has a cocompact lattice
and if $(\flow t)_t$ is completely conservative\footnote{That is,
  every wandering set has measure zero.} for the Liouville measure,
then
\begin{center}
  \fcolorbox{red}{white}{ $\displaystyle m_{F_{\rm su}} = m_{\rm Liou}$. }
\end{center}
\etheo

The following result, due to Bowen and Ruelle when $M$ is compact and
to Otal-Peigné \cite{OtaPei04} when $F=0$, completely solves the
problems of existence, uniqueness and explicit construction of
equilibrium states, see \cite[\S 6]{PauPolSha15}.

\btheo[Paulin-Pollicott-Schapira] Assume that $X=\wt M$ has pinched
sectional curvatures with uniformly bounded derivatives.\footnote{This
  assumption on the derivatives was forgotten in the statements of
  \cite{OtaPei04, PauPolSha15}.} If the Gibbs measure $m_F$ is finite,
then $\overline{m_F}=\frac{m_F}{\|m_F\|}$ is the unique equilibrium
state. Otherwise, there is no equilibrium state.
\etheo

We refer to Section \ref{subsect:variational} for an analogous
statement when $X$ is a tree, whose proof uses completely different
techniques.

\section{Basic ergodic properties of Gibbs measures}
\label{sect:ergodic}

We refer to \cite[Chap.~3, 5, 8]{PauPolSha15} and
\cite[Chap.~4]{BroParPau19} for details and complements on this
section.

\subsection{The Gibbs property}
\label{sect:Gibbsproperty}

In this section, we justify the terminology of Gibbs measures used above.

For every $\ell\in \Ga\bs\G X$, say $\ell=\Ga\wt \ell$, for every
$r>0$ and for all $t,t'\geq 0$, the (Bowen or) {\it dynamical ball}
$B(\ell;t,t',r)$ in the phase space $\Ga\bs\G X$ centered at $\ell$
with parameters $t,t',r$ is the image in $\Ga\bs\G X$ of the set of
geodesic lines in $\G X$ following the lift $\wt \ell$ at distance
less than $r$ in the time interval $[-t',t]$, that is, the image in
$\Ga\bs\G X$ of
\begin{center}
  \fcolorbox{blue}{white}{
    $\displaystyle B(\wt \ell;t,t',r)=\big\{\ell\,'\in \G X
    \;:\; \sup_{s\,\in\,[-t',\,t]}\;
    d_{X}(\,\wt \ell(s),\ell\,'(s)\,) <r\big\}$. }
\end{center}
The following definition of the Gibbs property is well adapted to the
possible noncompactness of the phase space $\Ga\bs \G X$. A $(\flow
t)_t$-invariant measure $m'$ on $\Ga\bs \G X$ satisfies the {\it Gibbs
  property} for the potential $F$ with {\it Gibbs constant} $c(F)\in
\RR$ if for every compact subset $K$ of $\Ga\bs \G X$, there exists
$r>0$ and $c_{K,r}\geq 1$ such that for all $t,t'\geq 0$ large enough,
for every $\ell$ in $\Ga\bs \G X$ with $\flow{-t'}\ell, \flow{t}
\ell\in K$, we have
\begin{center}
  \fcolorbox{green}{white}{
    $\displaystyle \frac{1}{C_{K,r}}\leq
    \frac{m'\big(B(\ell;t,t',r)\big)}
         {e^{\int_{-t'}^t \left(\,F(\flow{t}\ell)-c(F)\,\right)\,dt}}
         \leq C_{K,r}$. }
\end{center}
The following result is due to \cite[\S 3.8]{PauPolSha15} when $X=\wt
M$ and \cite[\S 4.2]{BroParPau19} in general.

\bprop\label{prop:gibbsgibbs}
The Gibbs measure $m_F$ satisfies the Gibbs property for $F$
with Gibbs constant $c(F)$ equal to the critical exponent $\delta_F$.
\eprop

Let us give a sketch of its proof, which explains the decorrelation of
the influence of the two points at infinity of the geodesic lines,
using the fact that the Gibbs measure is absolutely continuous with
respect to a product measure in the Hopf parametrisation. The key
geometric lemma is the following one.

\blemm For every $r>0$, there exists $t_r>0$ such that for all
$t,t'\geq t_r$ and $\ell\in\G X$, we have, using the Hopf
parametrisation at the footpoint $\ell(0)$,
$$
\OOO_{\ell(0)}B(\ell(-t'),r)\times \OOO_{\ell(0)}B(\ell(t),r)\times
\;]-1,1[\;\; \subset\; B(\ell;t,t',2r+2)
$$
$$
B(\ell;t,t',r)\;\subset \;\OOO_{\ell(0)}B(\ell(-t'),2r)\times
\OOO_{\ell(0)}B(\ell(t),2r)\times \;]-r,r[\;.
$$
\elemm

Let us give a proof-by-picture of the first claim, the second one
being similar. See the following picture. If a geodesic line $\ell'$
has its points at infinity $\ell'_-$ and $\ell'_+$ in the shadows seen
from $\ell(0)$ of $B(\ell(-t'),r)$ and $B(\ell(-t'),r)$ respectively,
then by the properties of triangles in negatively curved spaces, if
$t$ and $t'$ are large, then the image of $\ell'$ is close to the
union of the images of the geodesic rays from $\ell(0)$ to $\ell_-$
and $\ell_+$. The control on the time parameter in Hopf
parametrisation then says that $\ell'$ is staying at bounded distance
from $\ell$ in the time interval $[-t',t]$.

\begin{center}
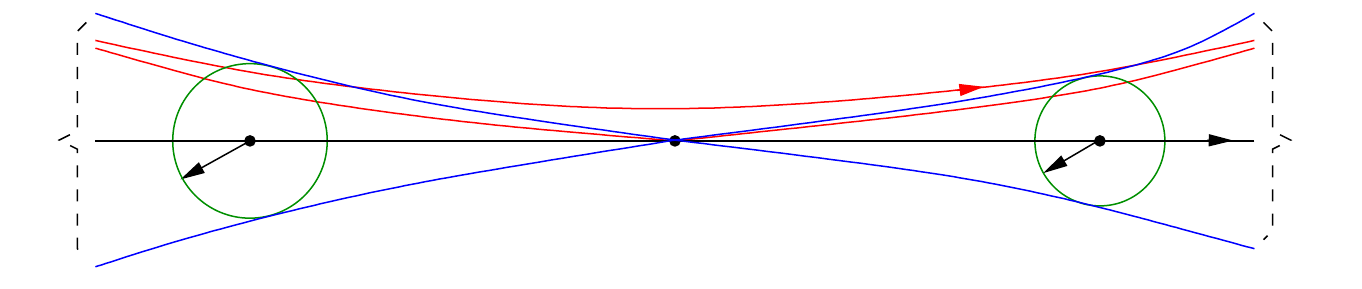
\end{center}

We now conclude the proof of Proposition \ref{prop:gibbsgibbs} by
using the boundedness of the Gibbs cocycles $C$ and $C^\iota$ on a
given compact subset $K$ in order to control the denominator in the
formula giving $\wt m_F$, and by using Mohsen's shadow lemma (see
Equation (2)) which estimates the Patterson measures of shadows of
balls.

\subsection{Ergodicity}
\label{sect:ergodicity}

In this section, we study the ergodicity property of the Gibbs
measures under the geodesic flow in the phase space.

The {\it Poincaré series} of the potential $F$ is 
\begin{center}
\fcolorbox{blue}{white}{
  $\displaystyle
  Q_F(s)= \sum_{\ga\in\Ga}\;\;\exp\Big(\int_{x_*}^{\ga x_*} (\wt F-s)\Big)\;.
  $ }
\end{center}
It depends on the basepoint $x_*$, but its convergence or divergence
does not. It converges if $s>\delta_F$ and diverges for $s<\delta_F$,
by the definition of the critical exponent $\delta_F$.

The following result has a long history, and we refer for instance to
\cite[\S 5]{PauPolSha15} and \cite[\S 4.2]{BroParPau19} for proofs,
and proofs of its following two corollaries.

\btheo[Hopf-Tsuji-Sullivan-Roblin] The following assertions are
equivalent.
\begin{enumerate}
\item The Poincaré series of $F$ diverges at the critical exponent of
  $F$~: $Q_F(\delta_F)=+\infty$.
\item The group action $(\partial_\infty X\times \partial_\infty
  X-{\rm Diag},\mu^\iota_{x_*}\otimes \mu_{x_*},\Ga)$ is ergodic and
  completely conservative.
\item The geodesic flow on the phase space with the Gibbs measure
  $(\Ga\bs\G X,m_F,(\flow t)_t))$ is ergodic and completely
  conservative.
\end{enumerate}
\etheo

\bcoro If $Q_F(\delta_F)=+\infty$, then there exists a Patterson
density for $F$, unique up to a positive scalar. It is atomless, and
the diagonal in $\partial_\infty X\times \partial_\infty X$ has
measure $0$ for the product measure $\mu^\iota_{x_*}\otimes
\mu_{x_*}$.
\ecoro

Let us give a sketch of the very classical proof of the first claim of
this corollary.

\medskip
\noindent{\bf Existence. } Using the properties of negatively curved
spaces, one can prove, denoting by $\D_x$ the Dirac mass at a point
$x$, that one can take
\begin{center}
  \fcolorbox{green}{white}{ $\displaystyle
    \mu_x=\lim_{s_i\ra\,\delta_F^+}\;\frac{1}{Q_F(s_i)}\;\sum_{\ga\in\Ga}
    \;\;\exp\Big(\int_{x}^{\ga x_*} (\wt F-s_i)\Big)\;\;\D_{\ga x_*}$, }
\end{center}
where the atomic measure before taking the limit is, when $x=x_*$, a
probability measure, hence has, for some sequence $(s_i)_{i\in\NN}$ in
$]\delta_F,+\infty[$ converging to $\delta_F$, a weakstar converging
    subsequence in the compact space of probability measures on the
    compact space $X\cup\partial_\infty X$.

\medskip
\noindent{\bf Uniqueness. } Let $(\mu'_x)_{x}$ be another Patterson
density. Up to positive scalars, we may assume that $\mu_{x_*}$ and
$\mu'_{x_*}$ are probability measures. Then $(\omega_x=\frac{1}{2}
(\mu_x+\mu'_x))_{x}$ is a Patterson density, $\mu_{x_*}$ is absolutely
continuous with respect to $\omega_{x_*}$, and by ergodicity, the
Radon-Nikodym derivative $\frac{d\mu_{x_*}}{d\,\omega_{x_*}}$ is
almost everywhere constant, hence the probability measures $\mu_{x_*}$
and $\omega_{x_*}$ are equal, hence $\mu_{x_*}=\mu'_{x_*}$.

\bcoro If $m_F$ is finite, then $Q_F(\delta_F)=+\infty$ (hence $(\flow
t)_t)$ is ergodic) and the normalised Gibbs measure $\overline{m_F}=
\frac{m_F}{\|m_F\|}$ is a cohomological invariant of the potential $F$.
\ecoro

\subsection{Mixing}
\label{sect:mixing}

In this section, we study the mixing property of the Gibbs measures
under the geodesic flow in the phase space. Recall that the {\it
  length spectrum} for the action of $\Ga$ on $X$ is the subgroup of
$\RR$ (hence of $\ZZ$ when $X$ is a tree) generated by the set of
lengths of the closed geodesic in $\Ga\bs X$ (or, in dynamical terms,
of the set of lengths of periodic orbits of the geodesic flow on the
phase space). See for instance \cite[\S 8.1]{PauPolSha15} when $X=\wt
M$ and \cite[\S 4.4]{BroParPau19} when $X$ is a tree for a proof of
the following result, which crucially uses the fact that the Gibbs
measure is absolutely continuous with respect to a product measure in
the Hopf parametrisation.

\btheo[Babillot] If the Gibbs measure $m_F$ is finite, then the
following assertions are equivalent.
\begin{enumerate}
\item The Gibbs measure $m_F$ is mixing under the geodesic flow
  $(\flow t)_t$.
\item The geodesic flow $(\flow t)_t$ is topologically mixing on its
  nonwandering set in the phase space.
\item The length spectrum of $\Ga$ is dense in $\RR$ if $X=\wt M$ or
  equal to $\ZZ$ if $X$ is a tree.
\end{enumerate}
\etheo

We summarise in the following result the known properties of the rate
of mixing of the geodesic flow in the manifold case when $X=\wt M$
(see \cite[\S 9.1]{BroParPau19}), refering to Section
\ref{sect:mixingrate} for the tree case, whose proof turns out to be
quite different.

Let $\alpha\in\;]0,1]$ and let $\C_{\rm b}^\alpha (Z)$ be the Banach
space\footnote{Recall that its norm (taking into account the
possible noncompactness of $Z$) is given by
$$
  \|f\|_\alpha= \|f\|_\infty +
  \sup_{\substack{x,\,y\,\in \,Z\\ 0<d(x,\,y)\leq 1}}
  \frac{|f(x)-f(y)|}{d(x,y)^\alpha}\;.
$$} of bounded $\alpha$-H\"older-continuous functions on a metric
space $Z$. When $X=\wt M$, we will say that the (continuous time)
geodesic flow on the phase space $T^1M=\Ga\backslash T^1 \wt M$ is
{\em exponentially mixing for the $\alpha$-H\"older regularity} or
that it has {\em exponential decay of $\alpha$-H\"older correlations}
for the potential $F$ if there exist $c',\kappa >0$ such that for all
$\phi,\psi\in \C_{\rm b}^\alpha (T^1 M)$ and $t\in\RR$, we have
$$
  \Big|\int_{T^1 M}
  \phi\circ\flow{-t}\;\psi\;d\overline{m_{F}}-
\int_{T^1 M}\phi\; d\overline{m_{F}}
\int_{T^1 M}\psi\;d\overline{m_{F}}\;\Big|
\le c'\;e^{-\kappa|t|}\;\|\phi\|_\alpha\;\|\psi\|_\alpha\,.
$$

\btheo Assume that $X=\wt M$ and that $M=\Ga\bs \wt M$ is
compact. Then the geodesic flow on the phase space $T^1M$ has
exponential decay of H\"older correlations if

\begin{itemize}
\item $M$ is two-dimensional, by \cite{Dolgopyat98},

\item $M$ is $1/9$-pinched and $F=0$, by
\cite[Coro.~2.7]{GiuLivPol13},

\item the potential $F$ is the unstable Jacobian $F^{\rm su}$, so
  that, up to a positive scalar, $m_F$ is the Liouville measure
  $m_{\rm Liou}$, by \cite{Liverani04}, see also \cite{Tsujii10},
  \cite[Coro.~5]{NonZwo15} who give more precise estimates,
 
\item $M$ is locally symmetric by \cite{Stoyanov11}, see also
  \cite{MohOh15} for some noncompact cases.
\end{itemize}
\etheo

Note that this gives only a very partial picture of the rate of mixing
of the geodesic flow in negative curvature, and it would be
interesting to have a complete result. Stronger results exist for the
Sobolev regularity when $\wt M$ is a symmetric space, $F=0$ and $\Ga$
is an arithmetic lattice (the Gibbs measure then coincides, up to a
multiplicative constant, with the Liouville measure): see for instance
\cite[Theorem~2.4.5]{KleMar96}, using spectral gap properties given by
\cite[Theorem 3.1]{Clozel03}. But this still does not give a complete
answer.

\section{Coding and rate of mixing for geodesic flows on trees}
\label{sect:mixingrate}

We refer to \cite[Chap.~5 and 9.2]{BroParPau19} for details and
complements on this section.

From now on, we assume that $X$ is (the geometric realisation of) a
simplicial tree $\XX$, and we write $\G\XX$ instead of $\G X$. We
consider the discrete group $\Ga$, the system of conductances $\wt c$
and the associated potential $F$ on the phase space $\Ga\bs \G\XX$ as
introduced in Section \ref{sect:construction}.

The study of the rate of mixing of the (discrete time) geodesic flow
on the phase space uses coding theory. But since, as explained, we make
no assumption of compactness on the phase space, and no hypothesis of
being without torsion on the group $\Ga$ in the huge class of examples
described in Section \ref{sect:construction}, the coding theory
requires more sophisticated tools than subshifts of finite type.

\subsection{Coding}
\label{subsect:coding}

Let $\A$ be a countable discrete set, called an {\it alphabet}, and
let $A= (A_{i,\,j})_{i,\,j\in\A}$ be an element in
$\{0,1\}^{\A\times\A}$, called a {\it transition matrix}. The
(two-sided, countable state) {\it topological shift}\footnote{We
  prefer not to use the frequent terminology of {\it topological
    Markov shift} as it could be misleading, many probability measures
  invariant under general topological shifts do not satisfy the Markov
  chain property that the probability to pass from one state to
  another depends only on the previous state, not of all past states.}
with alphabet $\A$ and transition matrix $A$ is the topological
dynamical system $(\Sigma,\sigma)$, where $\Sigma$, called the {\it
  shift space}, is the closed subset of the topological product space
$\A^\ZZ$ of {\it $A$-admissible} two-sided infinite sequences, defined
by
$$
\Sigma=\big\{x=(x_n)_{n\in\ZZ}\in \A^\ZZ\;:\; \forall \;n\in\ZZ,\;\;\;
A_{x_n,x_{n+1}}=1\}\;,
$$ 
and $\sigma: \Sigma\ra \Sigma$ is the (two-sided) {\it
  shift}\index{shift} defined by
$$
\forall\;x\in \Sigma,\;\forall\;n\in\ZZ,\;\;\;\;(\sigma(x))_n=x_{n+1}\;.
$$ 
We endow $\Sigma$ with the
distance
$$
d(x,x')= \exp\big(-\sup\big\{n\in\NN\;:\;\;
\forall\, i\,\in\,\{-n,\dots,n\},\;\;x_i\;=\;x'_i\big\}\,\big)\;.
$$

Let us denote by $\YY$ the (countable) quotient graph\footnote{The
  fact that the canonical projection is a morphism of graphs is the
  reason why we assumed $\Ga$ to be acting without mapping an edge to
  its inverse.}  $\Ga\bs \XX$. For every vertex or edge $x\in V\YY\cup
E\YY$, we fix a lift $\wt x$ in $V\XX\cup E\XX$, and we define
$G_x=\Ga_{\wt x}$ to be the stabiliser of $\wt x$ in $\Ga$.

\medskip
\noindent
\begin{minipage}{8.9cm}
For every $e\in E\YY$, we assume that $\widetilde{\overline{e}}
=\overline{\wt{e}}$. But there is no reason in general that
$\widetilde{t(e)} =t(\wt{e}\,)$. We fix $g_e\in\Ga$ mapping
$\widetilde{t(e)}$ to $t(\wt e\,)$ (which does exist), and we denote
by $\rho_e: G_e= \Ga_{\widetilde{e}}\ra \Ga_{\widetilde{t(e})} =
G_{t(e)}$ the conjugation $g\mapsto g_e^{-1}\,g\,g_e$ by $g_e$ on
$G_e$ (noticing that the stabiliser $\Ga_{\widetilde{e}}$ is contained in
the stabiliser $\Ga_{t(\widetilde{e})}$).
\end{minipage}
\begin{minipage}{6cm}
\begin{center}
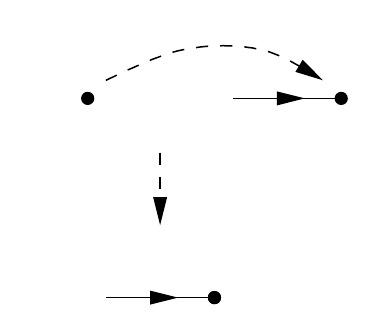
\end{center}
\end{minipage}

Let us try to code a geodesic line in the phase space $\Ga\bs\G\XX$.
The natural starting point is to write it as $\Ga\ell$ for some
$\ell\in\G X$, that is, to choose one of its lifts. We then have to
construct a coding which is independent of the choice of this
lift. For every $i\in\ZZ$, let us denote by $f_i=\ell([i,i+1])$ the
$i$-th edge followed by $\ell$, and by $e_i$ (also denoted by
$e_{i+1}^-(\ell)$ for later use) its image by the canonical $p:\XX\ra
\YY=\Ga\bs\XX$, which seems fit to be a natural part of the coding of
$\ell$. Since we will need to translate through our coding the fact
that $\ell$ is geodesic, hence has no backtracking, the edge $e_{i+1}$
(also denoted by $e_{i+1}^+(\ell)$ for later use) following $e_i$
seems to have a role to play.

\medskip
\begin{center}
  \!\!\!\!\!\!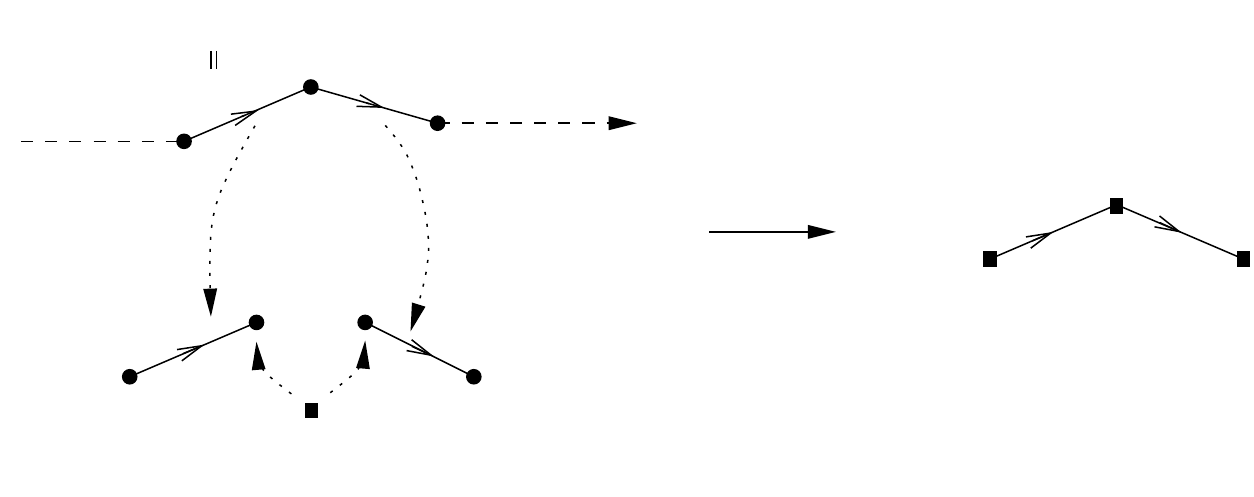
\end{center}
Since the terminal point of $f_i$ is the original point of $f_{i+1}$,
the terminal point of $e_i$ is naturally also the original point of
$e_{i+1}$.  But there is no reason for the terminal point of the
choosen lift $\widetilde{e_i}$ to also be the original point of the
choosen lift $\widetilde{e_{i+1}}$. Since $f_i$ and $\widetilde{e_i}$
both map by $p$ to $e_i$, we may fix $\ga_i\in\Ga$ such that $\ga_i
f_i=\widetilde{e_i}$, for every $i\in\ZZ$.

Now, note that the vertex stabilizers in $\Ga$ of vertices of $\XX$
are in general nontrivial (and we explained in Section
\ref{sect:construction} that it is important to allow them to become
very large in order to have numerous dynamically interesting
noncompact quotients of simplicial trees). The construction (see the
above diagram) provides a natural element $g_{e_i}^{\;-1}\,\ga_i\,
\ga_{i+1}^{\;-1} \,g_{\;\overline{e_{i+1}}}$ which stabilises the
lifted vertex $\widetilde{t(e_i)}$, hence belongs to $G_{t(e_i)}$.
Since we made choices for the elements $\ga_i$, the element
$g_{e_i}^{\;-1}\,\ga_i\, \ga_{i+1}^{\;-1} \,g_{\;\overline{e_{i+1}}}$
gives a well-defined double class $h_{i+1}(\ell)$ in
$\rho_{e_i}(G_{e_i})\bs G_{t(e_i)}/\rho_{\,\overline{e_{i+1}}}
(G_{e_{i+1}})$, which also seems fit to be another natural piece of
the coding of $\ell$.

It turns out that this construction is indeed working.  We take as
alphabet the (countable) set
$$
\A=\Big\{(e^-,h,e^+)\;:\;\begin{array}{l}
e^\pm\in E\YY {\rm ~with~} t(e^-)=o(e^+)\\
h\in \rho_{e^-}(G_{e^-})\bs G_{o(e^+)}/\rho_{\,\overline{e^+}}(G_{e^+})
{\rm ~with~} h\neq [1] {\rm ~if~} \overline{e^+}=e^-
\end{array}\Big\}\;.
$$ This last assumption of conditional nontriviality of the double
class codes the fact that $\ell$ being a geodesic line, the edge
$f_{i+1}$ is not the opposite edge of $f_i$, though $e_{i+1}$ might be
the opposite edge of $e_i$. And since in the tree $\XX$, being locally
geodesic implies being geodesic, it is very reasonable that we have
captured through our coding all the geodesic properties of the
geodesic lines and translated them into symbolic terms.  We take as
transition matrix over the alphabet $\A$ the matrix with entries
$$
A_{(e^-,\,h,\,e^+),\,({e'}^-,\,h',\,{e'}^+)} =\left\{\begin{array}{l}
1 {\rm ~if~}  e^+={e'}^-\\0 {\rm ~otherwise},\end{array}\right.
$$
which just says that we are glueing together the coding of pairs of
consecutive edges of the geodesic line. Note that since the tree is
locally finite, the transition matrix $A$ has finitely many nonzero
entries on each row and column, hence the associated shift space
$\Sigma$ is locally compact.

We then refer to \cite[\S 5.2]{BroParPau19} for a proof of the
following result, though almost everything is in the above picture!
We denote by $F_{\rm symb}:\Sigma\ra\RR$ the locally constant map
which associates to $\big((e^-_i, h_i, e^+_i)\big)_{i\in\ZZ}$ the
image $\wt c(\widetilde{e^+_0})$ by the system of conductances of the
lift of its first edge.

\btheo\label{theo:coding} The map
$$
\Theta:\left\{\begin{array}{ccl}
\Ga\bs\G\XX & \longrightarrow & \Sigma\\
\Ga\ell & \mapsto &
\big((e^-_i(\ell), h_i(\ell), e^+_i(\ell))\big)_{i\in\ZZ}
\end{array}\right.
$$
is a bilipschitz homeomorphism, conjugating the time $1$ map of
the (discrete time) geodesic flow $(\flow t)_{t\in\ZZ}$ to the shift
$\sigma$. Furthermore,
\begin{enumerate}
\item $(\Sigma,\sigma)$ is topologically transitive,\footnote{This
  comes from the assumption that there is no nontrivial proper
  $\Ga$-invariant subtree in $\XX$, since then $\partial_\infty X=
  \Lambda\Ga$, implying that the nonwandering set of the geodesic flow
  $(\flow t)_{t\in\ZZ}$ is the full phase space $\Ga\bs\G\XX$.}
\item if the Gibbs measure $m_F$ is finite and if the length spectrum
  of $\Ga$ is equal to $\ZZ$, then the probability measure
  $\PP=\Theta_*\overline{m_F}$ is mixing for the shift $\sigma$ on
  $\Sigma$,
\item the measure $\PP$ satisfies the Gibbs property on
  $(\Sigma,\sigma)$ with Gibbs constant $\delta_F$ for the potential
  $F_{\rm symb}$.\footnote{That is, with a formulation adapted to the
  possibility that the alphabet $\A$ may be infinite, for every finite
  subset $E$ of the alphabet $\A$, there exists $C_E\geq 1$ such that
  for all $p\leq q$ in $\ZZ$ and for every
  $x=(x_n)_{n\in\ZZ}\in\Sigma$ such that $x_p,x_q\in E$, we have
$$
\frac {1}{C_E}\le\frac{\PP([x_{p}, x_{p+1},\dots,x_{q-1}, x_{q}])}
{e^{-\delta_F(q-p+1)+\sum_{n=p}^{q}F_{\rm symb}(\sigma^n x)}} \le C_E\;.
$$
where $[x_{p}, x_{p+1},\dots,x_{q-1}, x_{q}]$ is the cylinder
$\{(y_n)_{n\in\ZZ}\in \Sigma\;:\;{\rm if~} p\leq n\leq q
{\rm ~then~} y_n=x_n\}$.}
\item if $(Z_n:x\mapsto x_n)_{n\in\ZZ}$ is the canonical random
  process in symbolic dynamics, then the pair $((Z_n)_{n\in\ZZ},\PP)$
  is not always a Markov chain.
  \end{enumerate}
\etheo

This last claim has lead to an erratum in the paper \cite{Kwon15}. The
pair $((Z_n)_{n\in\ZZ},\PP)$ is not a Markov chain for instance in
Example (2) at the beginning of Section \ref{sect:construction}, when
$\XX=\XX_q$ and $\Ga=\PGL_2(\FF_q[Y])$.\footnote{As noticed by
  J.-P.~Serre \cite{Serre83}, the image of almost every geodesic line
  of $\XX$ in the quotient ray $\Ga\bs X$ is a broken line which makes
  infinitely many back-and-forths from the origin of the quotient ray.

\begin{center}
\begin{picture}(0,0)%
\includegraphics{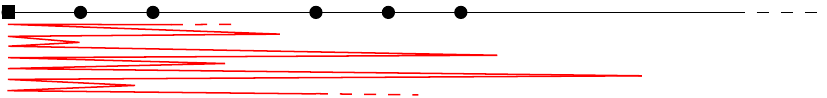}%
\end{picture}%
\setlength{\unitlength}{3812sp}%
\begingroup\makeatletter\ifx\SetFigFont\undefined%
\gdef\SetFigFont#1#2#3#4#5{%
  \reset@font\fontsize{#1}{#2pt}%
  \fontfamily{#3}\fontseries{#4}\fontshape{#5}%
  \selectfont}%
\fi\endgroup%
\begin{picture}(4102,466)(1311,-1743)
\end{picture}%

\end{center}
\nopagebreak
    There is absolutely no way to predict the probability of behaviour
    of the geodesic line image at a given time in terms of its recent
    past probabilities (except that when it starts to go down, it has
    to go down all the way to the origin).}

\subsection{Variational principle for simplicial trees}
\label{subsect:variational}

The first corollary of the coding results in the previous section is
the following existence and uniqueness result of equilibrium states
for the geodesic flow on the phase space $\Ga\bs \G\XX$ for the
potential $F$. 

\bcoro\label{coro:varprinctree} If $m_F$ is finite, then
$\overline{m_F}=\frac{m_F}{\|m_F\|}$ is the unique equilibrium state
for $F$ under the geodesic flow $(\flow t)_{t\in\ZZ}$ on $\Ga\bs
\G\XX$, and furthermore
\begin{center}
  \fcolorbox{red}{white}{ $\displaystyle P_F=\delta_F$. }
\end{center}
\ecoro

We only give a sketch of a proof, refering to \cite[\S
  5.4]{BroParPau19} for a complete one. We use the coding given in
Theorem \ref{theo:coding} with its properties (in particular the fact
that it satisfies the Gibbs property for a symbolic potential related
to the potential $F$).

Let $(\Sigma,\sigma)$ be a topological shift, with countable alphabet
$\A$. A $\sigma$-invariant probability measure $m$ on $\Sigma$ is a
{\it weak\footnote{The terminology comes from the fact that the
    assumptions bear only on the periodic points of $\sigma$.}  Gibbs
  measure} for a map $\phi:\Sigma\ra \RR$ with Gibbs constant $c(m)
\in\RR$ if for every $a\in \A$, there exists a constant $c_a\geq 1$
such that for all $n\in\NN-\{0\}$ and $x$ in the cylinder $[a]=
\{y=(y_n)_{n\in\ZZ}\in\Sigma\;:\;y_0=a\}$ such that $\sigma^n(x) =x$,
we have
$$
\frac {1}{c_a}\le\frac{m([x_{0}, x_{1},\dots, x_{n-1}])}
{e^{\sum_{i=0}^{n-1}\;(\,\phi(\sigma^i x)-c(m)\,)}} \le c_a\;.
$$

The following result of Buzzi is proved in
\cite[Appendix]{BroParPau19}, with a much weaker regularity assumption
on $\phi$, and it concludes the proof of Corollary
\ref{coro:varprinctree}.

\btheo[Buzzi] Let $(\Sigma,\sigma)$ be a topological shift and
$\phi:\Sigma\ra\RR$ a bounded Hölder-continuous function. If $m$ is a
weak Gibbs measure for $\phi$ with Gibbs constant $c(m)$, then
$P_\phi=c(m)$ and $m$ is the unique equilibrium state for the
potential $\phi$.
\etheo

\subsection{Rate of mixing for simplicial trees}
\label{subsect:mixingrate}

Let us first recall the definition of an exponential mixing rate for
discrete time dynamical systems.

Let $(Z,m,T)$ be a dynamical system with $(Z,m)$ a metric probability 
space and let $T:Z\ra Z$ be a (not necessarily invertible) measure
preserving map. For all $n\in\NN$ and $\phi,\psi\in\LL^2(m)$, the
(well-defined) $n$-th {\it correlation coefficient} of $\phi,\psi$ is
$$
\operatorname{cov}_{m,\,n}(\phi,\psi)=
\int_{Z}(\phi\circ T^n)\;\psi\;dm-\int_{Z}\phi\; dm\;\int_{Z}\psi\;dm\;.
$$
Let $\alpha\in\;]0,1]$. As for the case of flows in Section
\ref{sect:mixing}, we will say that the dynamical system $(Z,m,T)$ is
{\it exponentially mixing for the $\alpha$-H\"older regularity} or
that it has {\it exponential decay of $\alpha$-H\"older
correlations} if there exist $c',\kappa >0$ such that for all
$\phi,\psi\in \C_{\rm b} ^\alpha(Z)$ and $n\in\NN$, we have
$$
|\operatorname{cov}_{m,\,n}(\phi,\psi)|
\le c'\;e^{-\kappa\, n }\;\|\phi\|_\alpha\;\|\psi\|_\alpha\,.
$$
Note that this property is invariant under measure preserving
conjugations of dynamical systems by bilipschitz homeomorphisms. In
our case, $T$ will be either the time $1$ map of the geodesic flow
$(\flow t)_{t\in\ZZ}$ on the phase space $Z=\Ga\bs\G\XX$ or the
two-sided shift $\sigma$ on a two-sided topological shift space
$\Sigma$ or (see below) the one-sided shift $\sigma_+$ on a one-sided
topological shift space $\Sigma_+$.

\medskip
The following result is one of the new results contained in the book
\cite{BroParPau19}. For every finite subset $E$ in $\Ga\bs V\XX$, let
$\tau_E:\Ga\bs\G\XX\ra \NN\cup\{+\infty\}$ be the first positive
passage time of geodesic lines in $E$, that is, the map
$$
\ell\mapsto\inf\{n\in\NN-\{0\}\;:\; \flow{n}\ell(0)\in E\}\;.
$$
The following result says that if the tree quotient contains a
finite subset in which the geodesic lines with large return times have
an exponentially decreasing mass, then the (discrete time) geodesic
flow on the phase space has exponential decay of correlations. This
condition turns out to be quite easy to check on practical examples,
see for instance \cite[\S 9.2]{BroParPau19}. 

\btheo\label{theo:mixingratetree}
If $m_F$ is finite and mixing for $(\flow t)_{t\in\ZZ}$, if
there exist a finite subset $E$ in $\Ga\bs V\XX$ and $c'',\kappa'>0$
such that
$$
\forall\;n\in\NN,\;\;\;m_F(\{\ell\in\Ga\bs\G\XX\;:\;
\ell(0)\in E, \tau_E(\ell)\geq n\})\leq c''e^{-\kappa' n}\;,
$$
then for every $\alpha\in\;]0,1]$, the (discrete time) dynamical system
$(\Ga\bs\G\XX,m_F,(\flow t)_{t\in\ZZ})$ is exponentially mixing for
the $\alpha$-Hölder regularity.
\etheo

The hypothesis of Theorem \ref{theo:mixingratetree} is for instance
satisfied for Example (2) at the beginning of Section
\ref{sect:construction} with $\XX=\XX_q$ and $\Ga =\PGL_2(\FF_q[Y])$,
taking $E$ consisting of the origin of the modular ray $\Ga\bs \XX_q$,
and using the exponential decay of the stabilisers orders along a lift
of the modular ray in $\XX_p$. In this case, the quotient graph
$\Ga\bs\XX$ has linear growth. We gave in \cite[page 193]{BroParPau19}
examples where the quotient graph $\Ga\bs\XX$ has exponential
growth.

Here is an example where the quotient graph has quadratic growth, for
every even $q\geq 2$. The tree $\XX$ is the regular tree of degrees
$q+2$. The vertex group of the top-left vertex $x_*$ of the quotient
graph is $\ZZ/(\frac{q}{2}+1)\ZZ$. A set $E$ as in Theorem
\ref{theo:mixingratetree} consists of the three vertices at distance
at most $1$ from $x_*$. The vertex group of a vertex at distance
$k\geq 1$ from $x_*$ on the left vertical ray is $\ZZ/(q+1)^k\ZZ$. The
vertex group of a vertex not on the left vertical ray, at distance
$k\geq 1$ from $x_*$ is $\ZZ/q\ZZ\times\ZZ/(q+1)^{k-1}\ZZ$. The number
at the beginning of each edge represents the index of the edge group
inside the vertex group of its origin.

\begin{center}
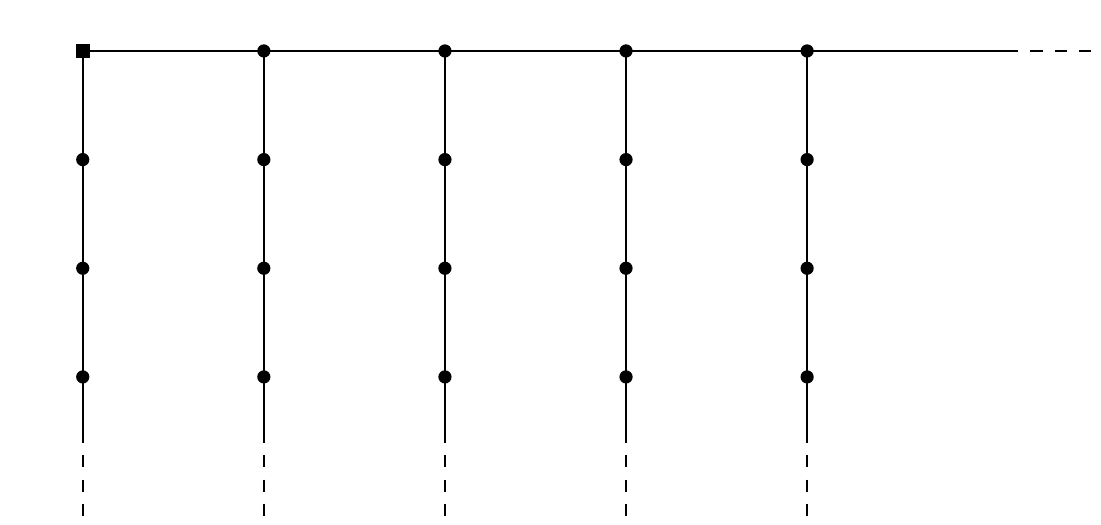
\end{center}

Recall that two {\it growth functions} $f$ and $f'$, that is, two
increasing maps from $\NN$ to $\NN-\{0\}$, are {\it equivalent} if
there exist two integers $c\geq 1$ and $c'\geq 0$ such that for every
$n\in\NN$ large enough, we have $f(\lfloor \frac{1}{c}\,n-c'\rfloor)
\leq f'(n)\leq f(c\,n+c')$.  The {\it type of growth} of an infinite,
connected, locally finite graph $Y$ is the equivalence class of the
map $n\mapsto \card \;B_{VY}(v_0,n)$, which does not depend on the
choice of a base point $v_0\in VY$, nor on the quasi-isometry type of
$Y$.
  
It is well known (see for instance \cite{Choucroun94b,Hughes04} or
\cite[\S 6.2]{GriNekSus00}) that every totally disconnected compact
metric space is homeomorphic to the boundary at infinity of a
simplicial tree with uniformly bounded degrees, and that any
increasing positive integer sequence $(a_n)_{n\in\NN}$ with at most
exponential speed (that is, there exists $k\in\NN$ such that
$a_{n+1}\leq ka_{n}$ for every $n\in\NN$) is, up to the above
equivalence, the sequence of orders of the balls of an infinite rooted
simplicial tree with uniformly bounded degrees. Hence the following
result (not contained in \cite{BroParPau19}) says that we can realize
any space of ends, or any at most exponential type of growth, in the
quotient graph of an action of a group on a tree satisfying the
hypothesis of Theorem \ref{theo:mixingratetree}.

\bprop\label{prop:allendsallgrowth}
For every rooted tree $(\T,*)$ with uniformly bounded degrees,
there exists a simplicial tree $\XX$ and a discrete group $\Ga$ of
automorphisms of $\XX$ as in the beginning of Section
\ref{sect:construction} such that $\Ga$ is a lattice, $\Ga\bs \XX=\T$
and the geodesic flow $(\flow t)_{t\in\ZZ})$ is exponentially mixing
for the $\alpha$-Hölder regularity on $\Ga\bs \G\XX$ for the zero
potential.
\eprop

\dem We refer for instance to \cite[\S I.5]{Serre83} for
background on graphs of groups.

Let us fix $q\in\NN$ large enough compared with the maximum degree $d$
of $\T$. We define a graph of groups $(\T,G_\bullet)$ with underlying
graph $\T$ as follows. For every vertex $v$ of $\T$ at distance $n$ of
the root $*$, we define $G_v=\ZZ/q^{n+1}\ZZ$. For every edge $e$ whose
closest vertex to the root $*$ is at distance $n$ from $*$, we define
$G_e=\ZZ/q^{n+1}\ZZ$. For every edge $e$ pointing away from the root,
we define the monomorphism $G_e\ra G_{o(e)}$ to be the identity, and
the monomorphism $G_e\ra G_{t(e)}$ to be the identity on the first
factors, so that the index of $G_e$ in $G_{o(e)}$ is $1$ and the index
of $G_e$ in $G_{t(e)}$ is $q$.

Let $\Ga$ and $\XX$ be respectively the fundamental group (using the
root as the basepoint) and the Bass-Serre tree of the graph of groups
$(\T,G_\bullet)$. 
Then the degrees of the vertices of $\XX$ are at least $3$ (actually
at least $q$) and at most $q+d-1$, and for every $n$, we have
\begin{equation}\label{eq:expodecayreturn}
\sum_{x\in V\T\;:\; d(x,*)=n}\frac{1}{|G_x|}\leq d^n/q^n\;.
\end{equation}
Since $q$ is large compared to $d$, this implies that the volume of
$(\T,G_\bullet)$ is finite, hence $\Ga$ is a lattice.

Since the potential is the zero potential, the Gibbs measure is the
Bowen-Margulis measure, and up to a positive scalar, the Patterson
density is, by \cite[Prop.~4.16]{BroParPau19}, the Hausdorff measure
of the visual distance $d_{x_*}$. Since $q$ is large compared to $d$,
the set of points at infinity of lifts in $\XX$ of geodesic rays in
$\T$
starting from the root has measure $0$ for the Patterson
density. Since for every edge in $\T$ pointing away from the root, the
index of its edge group in its original vertex group is $1$, almost
every geodesic line for the Bowen-Margulis measure (which is
absolutely continuous with respect to the product measure of the
Patterson densities on its two endpoints and the counting measure
along its image) maps in $\T$ to a path making infinitely many
back-and-forth from the root.  If $E=\{*\}$ is the singleton in
$V\!\T$ consisting of the root, since $q$ is large compared to $d$,
Equation \eqref{eq:expodecayreturn} then shows that the hypothesis of
Theorem \ref{theo:mixingratetree} is satisfied, and this concludes the
proof of Proposition \ref{prop:allendsallgrowth}.  \cqfd

\medskip
We conclude this survey with a sketch of proof of Theorem
\ref{theo:mixingratetree}, sending to \cite[\S 9.2]{BroParPau19} for a
complete proof. We thank Omri Sarig for a key idea in the proof of
this theorem.

\medskip
\noindent{\bf Step 1. } The first step consists in passing from the
geometric dynamical system to a two-sided symbolic dynamical system,
using Section \ref{subsect:coding}.

\medskip
Let $\A,A,\Sigma,\sigma,\Theta,\PP$ be as given in Theorem
\ref{theo:coding} for the coding of the (discrete time) geodesic flow
on the phase space $\Ga\bs\G\XX$. Let $\pi_+:\Sigma\ra \A^\NN$ be the
natural projection defined by $(x_n)_{n\in\ZZ}\mapsto (x_n)_{n\in\NN}$,
let $(\Sigma_+,\sigma_+)$ be the one-sided topological shift
constructed as for the two-sided one with the same alphabet $\A$ and
same transition matrix $A$, with $\Sigma_+\subset \A^\NN$.  Let
$$
\E=\{(e^-,h,e^+)\in\A \;:\; t(e^-)=o(e^+)\in E\}
$$
which is a finite subset of the alphabet, and $\tau_{\E}:\Sigma_+ \ra
\NN$ the first positive passage time in $\E$ of the shift orbits, that
is, the map $x=(x_n)_{n\in\NN}\mapsto \inf\{n\in\NN-\{0\}\;:\;
x_n\in\E\}$.

The rate of mixing statement for two-sided symbolic dynamical system,
that we will prove in Step 2, is the following one.

\btheo \label{theo:critexpdecaysimpldynsymb} Let
$(\A,A,\Sigma,\sigma)$ be a locally compact transitive two-sided
topological shift, and let $\PP$ be a mixing $\sigma$-invariant
probability measure with full support on $\Sigma$. Assume that
\begin{enumerate}
\item[(1)] for every $n \in\NN$ and for every $A$-admissible finite
  sequence $w= (w_0,\dots,w_n)$ in $\A$, the (measure theoretic)
  Jacobian of the map
  $$
  f_w:\{(x_k)_{k\in\NN}\in\pi_+(\Sigma)\;:\;x_0=w_n\}\ra\{(y_k)_{k\in\NN}
  \in \pi_+(\Sigma)\;:\;y_0=w_0,\dots, y_n=w_n\}
  $$
  defined by $(x_0,x_1,x_2,\dots)\mapsto (w_0,\dots,w_n,
  x_1,x_2,\dots)$, with respect to the restrictions of the pushforward
  measure $(\pi_+)_*\PP$, is constant;
\item[(2)] there exist a finite subset $\E$ of $\A$ and $c'',\kappa'
  >0$ such that for every $n\in\NN$, we have
$$
\PP\big(\{x\in \Sigma\;:\;x_0\in \E \;{\rm and}\;
\tau_\E(x)\geq n\}\big)\leq c''\;e^{-\kappa' n}\;.
$$
\end{enumerate}
Then $(\Sigma,\sigma,\PP)$ has exponential decay of $\alpha$-H\"older
correlations.  
\etheo

Theorem \ref{theo:mixingratetree} follows from Theorem
\ref{theo:critexpdecaysimpldynsymb} by using the coding given in
Theorem \ref{theo:coding}. The verification of Assertion (2) is
immediate as it corresponds to the assumption of Theorem
\ref{theo:mixingratetree}. The one of Assertion (1) is a bit
technical, using a strengthened version of Mohsen's shadow lemma for
trees.

\medskip
\noindent{\bf Step 2. } The second step consists in passing from the
two-sided symbolic dynamical system to a one-sided symbolic dynamical
system.

Let $(\Sigma_+,\sigma_+)$ be the one-sided topological shift with the
same alphabet $\A$ and same transition matrix $A$ as the two-sided one
in the statement of Theorem \ref{theo:critexpdecaysimpldynsymb}, with
\mbox{$\Sigma_+=\pi_+(\Sigma)$,} and let $\PP_+=(\pi_+)_*\PP$. Recall
that the {\it cylinders} in $\Sigma_+$ are the subsets defined for
$k\in\NN$ and $w_0,\dots w_k\in\A$ by
$$
[w_0,\dots,w_k]=\{x=(x_n)_{n\in\NN}\in\Sigma_+\;:\; x_0=w_0,\dots,x_k=w_k \}\;.
$$

The rate of mixing statement for one-sided symbolic dynamical system,
that we will prove in Step 3, is the following one.

\btheo \label{theo:critexpdecaysimpldynsymbonsesided} Let
$(\A,A,\Sigma_+,\sigma_+)$ be a locally compact transitive one-sided
topological shift, and let $\PP_+$ be a mixing $\sigma$-invariant
probability measure with full support on $\Sigma_+$. Assume that
\begin{enumerate}
\item[(1)] for every $n\in\NN$ and for every $A$-admissible finite
  sequence $w= (w_0,\dots,w_n)$ in $\A$, the Jacobian of the map
  between cylinders
  $$
  f_w:[w_n]\ra[w_0,\dots,w_n]
  $$
  defined by $(x_0,x_1,x_2,\dots)\mapsto (w_0,\dots,w_n,
  x_1,x_2,\dots)$, with respect to the restrictions of $\PP_+$, is
  constant;
\item[(2)] there exist a finite subset $\E$ of $\A$ and $c'',\kappa'
  >0$ such that for every $n\in\NN$, we have
$$
\PP_+\big(\{x\in \Sigma_+\;:\;x_0\in \E \;{\rm and}\;
\tau_\E(x)\geq n\}\big)\leq c''\;e^{-\kappa' n}\;.
$$
\end{enumerate}
Then $(\Sigma_+,\sigma_+,\PP_+)$ has exponential decay of
$\alpha$-H\"older correlations.
\etheo

Theorem \ref{theo:critexpdecaysimpldynsymb} follows from Theorem
\ref{theo:critexpdecaysimpldynsymbonsesided} by a classical argument
due to Sinai and Bowen (and explained to the authors by Buzzi), saying
that if the one-sided symbolic dynamical system $(\Sigma_+,\sigma_+,
(\pi_+)_*\PP)$ is exponentially mixing, then so is the two-sided
symbolic dynamical system $(\Sigma,\sigma,\PP)$.

\medskip
\noindent{\bf Step 3. } The third and final step that we sketch is a
proof of Theorem \ref{theo:critexpdecaysimpldynsymbonsesided}, using
as main tool a Young's tower argument.

We implicitely throw away from $\Sigma_+$ the measure zero subset of
points $x\in\Sigma_+$ whose orbit under the shift $\sigma_+$ does not
pass infinitely many times in the open nonempty finite union of
fundamental cylinders
$$
\Delta_0=\bigcup_{a\in\E}\;\;[a]\;.
$$
We denote by $\Phi :\Sigma_+\ra\Delta_0$ the first positive time
passage map, defined by $x\mapsto \sigma_+^{\tau_\E(x)}(x)$. We denote
by $W$ the set of excursions outside $\E$, that is, the set of
$A$-admissible finite sequences $(w_0,\dots,w_n)$ in $\A$ such that
$w_0,w_n\in\E$ and $w_i\notin\E$ for $1\leq i\leq n-1$.

We have the following properties.
\begin{enumerate}
\item The set $\{[a]\;:\;a\in\E\}$ is a finite measurable partition of
  $\Delta_0$. For every $a\in\E$, the set $\{[w]\;:\; w\in W, w_0=a\}$
  is a countable measurable partition of $[a]$.
\item For every $w\in W$, the first positive passage time $\tau_\E$ is
  positive on every excursion cylinder $[w]$, and if $w_n$ is the last
  letter of $w$, then the restriction $\Phi:[w]\ra[w_n]$ is a
  bijection with constant with constant Jacobian with respect to
  $\PP_+$ (actually much less is needed in order to apply Young's
  arguments).
\item The first positive time passage map $\Phi$ satisfies strong
  dilations properties on the excursion cylinders. More precisely, for
  every excursion $w= (w_0,\dots,w_n)\in W$, for every $k\leq n-1$,
  for all $x,y\in [w]$, we have $d(\Phi(x),\Phi(y))\geq e\; d(x,y)$
  and $d(\sigma_+^kx,\sigma_+^ky))< d(\Phi(x),\Phi(y))$.
\end{enumerate}
Let us fix $\alpha\in\;]0,1]$. Then an adaptation of
\cite[Theo.~3]{Young99} implies that there exists $\kappa>0$ such
that for all $\phi,\psi\in \C_{\rm b} ^\alpha(\Sigma_+)$, there
exists $c_{\phi,\psi}>0$ such that for every $n\in\NN$, we have
$$
|\operatorname{cov}_{\PP_+,\,n}(\phi,\psi)|
\le c_{\phi,\psi}\;e^{-\kappa\, n }\;.
$$
An argument using the Principle of Uniform Boundedness due to Chazotte
then allows us to take $c_{\phi,\psi}=c'\;e^{-\kappa\, n }\;
\|\phi\|_\alpha\; \|\psi\|_\alpha$ for some constant $c'>0$.

{\small \bibliography{../biblio} }

\begin{thebibliography}{BAPP}

\bibitem[BaL]{BasLub01}
H.~Bass and A.~Lubotzky.
\newblock {\it Tree lattices}.
\newblock {Prog. in Math. {\bf 176}, Birkhäuser, 2001}.

\bibitem[BH]{BriHae99}
M.~R. Bridson and A.~Haefliger.
\newblock {\it Metric spaces of non-positive curva\-tu\-re}.
\newblock {Grund. math. Wiss. {\bf 319}, Springer Verlag, 1999}.

\bibitem[BPP]{BroParPau19}
A.~Broise-Alamichel, J.~Parkkonen and F.~Paulin.
\newblock {\it Equidistribution and counting under equilibrium states in
  negative curvature and trees. Applications to non-Archimedean Diophantine
  approximation}.
\newblock {With an Appendix by J.~Buzzi. Prog. Math. {\bf 329}, Birkhäuser,
  2019}.

\bibitem[Cho]{Choucroun94b} 
F.M.~Choucroun.
\newblock {\it Arbres, espaces ultramétriques et bases de structure uniforme}.
\newblock {Geom. Dedicata {\bf 53} (1994) 69--74}.

\bibitem[Clo]{Clozel03}
L.~Clozel.
\newblock {\it D\'emonstration de la conjecture $\tau$}.
\newblock {Invent. Math. {\bf 151} (2003) 297--328}.

\bibitem[Dol]{Dolgopyat98}
D.~Dolgopyat.
\newblock {\it On decay of correlation in Anosov flows}.
\newblock {Ann. of Math. {\bf 147} (1998) 357--390}.

\bibitem[GLP]{GiuLivPol13}
P.~Giulietti, C.~Liverani and M.~Pollicott.
\newblock {\it Anosov flows and dynamical zeta functions}.
\newblock {Ann. of Math. {\bf 178} (2013) 687--773}.

\bibitem[GNS]{GriNekSus00}
R.~Grigorchuk,  V.~Nekrashevich and V.~Sushchanskii. 
\newblock {\it Automata, dynamical systems, and groups}.
\newblock {Proc. Steklov Inst. Math. {\bf 231} (2000) 128--203}.

\bibitem[Hug]{Hughes04}
B.~Hughes.
\newblock {\it Trees and ultrametric spaces: a categorical equivalence}.
\newblock {Adv. Math. {\bf 189} (2004) 148--191}.

\bibitem[Kat]{Katok92}
S.~Katok.
\newblock {\it Fuchsian groups}.
\newblock {Univ. Chicago Press, 1992}.

\bibitem[KH]{KatHas95}
A.~Katok and B.~Hasselblatt.
\newblock {\it Introduction to the modern theory of dynamical systems}.
\newblock {Ency. Math. App. {\bf 54}, Camb. Univ. Press, 1995}.

\bibitem[KM]{KleMar96}
D.~Kleinbock and G.~Margulis.
\newblock {\it Bounded orbits of nonquasiunipotent flows on homogeneous
  spaces}.
\newblock {Sinai's Moscow Seminar on Dynamical Systems, 141--172, Amer. Math.
  Soc. Transl. Ser. {\bf 171}, Amer. Math. Soc. 1996}.

\bibitem[Kwo]{Kwon15}
S.~Kwon.
\newblock {\it Effective mixing and counting in Bruhat-Tits trees}.
\newblock {Erg. Theo. Dyn. Sys. {\bf 38} (2018) 257--283. Erratum {\bf 38}
  (2018) 284}.

\bibitem[Live]{Liverani04}
C.~Liverani.
\newblock {\it On contact Anosov flows}.
\newblock {Ann. of Math. {\bf 159} (2004) 1275--1312}.

\bibitem[Liv\v{s}]{Livsic72}
A.~Liv\v{s}ic.
\newblock {\it Cohomology of dynamical systems}.
\newblock {Math. USSR-Izv. {\bf 6} (1972) 1278--1301}.

\bibitem[Mar]{Margulis91}
G.~Margulis.
\newblock {\it Discrete subgroups of semi-simple groupes}.
\newblock {Ergeb. Math. Grenz. {\bf 17}, Springer Verlag, 1991}.

\bibitem[MO]{MohOh15}
A.~Mohammadi and H.~Oh.
\newblock {\it Matrix coefficients, counting and primes for orbits of
  geometrically finite groups}.
\newblock {J. Euro. Math. Soc. {\bf 17} (2015) 837--897}.

\bibitem[NZ]{NonZwo15}
S.~Nonnenmacher and M.~Zworski.
\newblock {\it Decay of correlations for normally hyperbolic trapping}.
\newblock {Invent. Math. {\bf 200} (2015) 345--438}.

\bibitem[OP]{OtaPei04}
J.-P. Otal and M.~Peign\'e.
\newblock {\it Principe variationnel et groupes kleiniens}.
\newblock {Duke Math. J. {\bf 125} (2004) 15--44}.

\bibitem[Pau]{Paulin97d}
F.~Paulin.
\newblock {\it On the critical exponent of discrete group of hyperbolic
  isometries}.
\newblock {Diff. Geom. and its App. {\bf 7} (1997) 231--236}.

\bibitem[PPS]{PauPolSha15}
F.~Paulin, M.~Pollicott, and B.~Schapira.
\newblock {\it Equilibrium states in negative curvature}.
\newblock {Astérisque {\bf 373}, Soc. Math. France, 2015}.

\bibitem[Rob]{Roblin03}
T.~Roblin.
\newblock {\it Ergodicit\'e et \'equidistribution en courbure n\'egative}.
\newblock {M\'emoire Soc. Math. France, {\bf 95} (2003)}.

\bibitem[Rue]{Ruelle04}
D.~Ruelle.
\newblock {\it Thermodynamic formalism: The mathematical structure of
  equilibrium statistical mechanics}.
\newblock {2nd ed., Cambridge Math. Lib. Cambridge Univ. Press, 2004}.

\bibitem[Ser]{Serre83}
J.-P. Serre.
\newblock {\it Arbres, amalgames, SL$_2$}.
\newblock {3\`eme \'ed. corr., Ast\'erisque {\bf 46}, Soc. Math. France, 1983}.

\bibitem[Sto]{Stoyanov11}
L.~Stoyanov.
\newblock {\it Spectra of Ruelle transfer operators for axiom A flows}.
\newblock {Nonlinearity {\bf 24} (2011) 1089--1120}.

\bibitem[Tsu]{Tsujii10}
M.~Tsujii.
\newblock {\it Quasi-compactness of transfer operators for contact Anosov
  flows}.
\newblock {Nonlinearity {\bf 23} (2010) 1495--1545}.

\bibitem[You]{Young99}
L.-S. Young.
\newblock {\it Recurrence times and rates of mixing}.
\newblock {Israel J. Math. {\bf 110} (1999) 153--188}.

\bibitem[Zem]{Zemanian91}
A.~Zemanian.
\newblock {\it Infinite electrical networks}.
\newblock {Cambridge Tracts Math. {\bf 101}, Cambridge Univ. Press, 1991}.

\end{thebibliography}

\bigskip
{\small

\noindent \begin{tabular}{l}
Laboratoire de math\'ematique d'Orsay, UMR 8628 CNRS\\
Universit\'e Paris-Saclay, Bâtiment 307,
91405 ORSAY Cedex, FRANCE\\
{\it e-mail: anne.broise@universite-paris-saclay.fr}
\end{tabular}

\medskip\noindent \begin{tabular}{l} 
Department of Mathematics and Statistics, P.O. Box 35\\ 
40014 University of Jyv\"askyl\"a, FINLAND.\\
{\it e-mail: jouni.t.parkkonen@jyu.fi}
\end{tabular}

\medskip\noindent \begin{tabular}{l}
Laboratoire de math\'ematique d'Orsay, UMR 8628 CNRS\\
Universit\'e Paris-Saclay, Bâtiment 307,
91405 ORSAY Cedex, FRANCE\\
{\it e-mail: frederic.paulin@universite-paris-saclay.fr}
\end{tabular}
}

\end{document}